\documentstyle[11pt]{article}

\newskip\stdskip                      
\newenvironment{prf}{\par{\bf Proof}\ }{\hfill$\Box$\par}
\newcommand{\cl}{\centerline}         
\newenvironment{changemargin}[2]{%
 \begin{list}{}{%
  \setlength{\topsep}{0pt}%
  \setlength{\leftmargin}{#1}%
  \setlength{\rightmargin}{#2}%
  \setlength{\listparindent}{\parindent}%
  \setlength{\itemindent}{\parindent}%
  \setlength{\parsep}{\parskip}%
 }%
\item[]}{\end{list}}

\def\SDP{\,  \vrule width0.4pt height 6pt depth 0pt \kern -4.5pt\times  } 

\def\sdp{\,    \vrule width0.4pt height 4pt depth 0pt \kern -1.5pt\times } 

\input xy
\xyoption{all}

\newtheorem{thrm}{Theorem}[section]
\newtheorem{lem}[thrm]{Lemma}
\newtheorem{cor}[thrm]{Corollary}
\newtheorem{rem}[thrm]{Remark}
\newtheorem{defn}[thrm]{Definition}
\newtheorem{exmpl}[thrm]{Example}
\newtheorem{prop}[thrm]{Proposition}

\def\maprt#1{\, \smash{\mathop{\longrightarrow}\limits^{#1}}\, }
\def\maplft#1{\smash{\mathop{\longleftarrow}\limits^{#1}}}

\def\bd{  \begin{diagram}    }
\def\ed{  \end{diagram}      }

\def\thm #1  {\medskip\noindent{\bf #1}\quad}
\def\pf #1 {\smallskip\noindent{\bf #1}\par\nobreak\noindent}
\def\term #1{{\bf #1}}

\font\tenmsb=msbm10 scaled \magstep 1
\font\sevenmsb=msbm7
\font\fivemsb=msbm5
\newfam\msbfam
\textfont\msbfam=\tenmsb
\scriptfont\msbfam=\sevenmsb
\scriptscriptfont\msbfam=\fivemsb

\def\Bbb#1{{\fam\msbfam\relax#1}}

\def\sec{\mathhexbox278}

\def\ZZ{{\Bbb Z}}

\def\CC{{\Bbb C}}

\def\s{\Sigma}
\def\smsh{\wedge}

\def\cross{\times}
\def\wdg{\vee}

\def\twdl{\widetilde}
\def\<{\langle}
\def\>{\rangle}

\def\sseq{\subseteq}

\def\bprp{\begin{prop}}
\def\eprp{\end{prop}}

\def\bthm{\begin{thrm}}
\def\ethm{\end{thrm}}

\def\blem{\begin{lem}}
\def\elem{\end{lem}}

\def\bcor{\begin{cor}}
\def\ecor{\end{cor}}

\def\brmk{\begin{rem}}
\def\ermk{\end{rem}}

\def\bdfn{\begin{defn}}
\def\edfn{\end{defn}}

\def\bexm{\begin{exmpl}}
\def\eexm{\end{exmpl}}

\def\|{\, \bigm|\, }

\def\A{{\cal A}}

\def\F{{\cal F}}

\def\S{{\cal S}}

\def\cp{\CC {\mathrm P} }

\def\kl{{\mathrm{kl}}}
\def\cl{{\mathrm{cl}}}
\def\cat{{\mathrm{cat}}}
\def\wcat{{\mathrm{wcat}}}
\def\kit{{\mathrm{kit}}}
\def\id{{\mathrm{id}}}

\def\l{\ell}

\def\L{{\cal L}}
\def\LL{L}

\def\dequivarrow{\ar@{-}@<.4ex>[d] \ar@{-}[d]  \ar@{-}@<-.4ex>[d]}
\def\requivarrow{\ar@{-}@<.4ex>[r] \ar@{-}[r]  \ar@{-}@<-.4ex>[r]}

\def\longdecom#1#2#3#4#5#6#7{
\xymatrix{
{#4}_0 \ar[d]   & & {#4}_{#6} \ar[d]   \\
{#2}_0 \ar@{=}[d] \ar[r] &\cdots\ar[r]& {#2}_{#6} \dequivarrow \ar[r]
    &\cdots\ar[r] &  {#2}_{{#7} + {#6}} \dequivarrow \\
{#1} \ar[rr] && {#2} \ar[rr]_{#5} && {#3} \\ }
}

\def\ldiag#1#2#3#4#5{  
\xymatrix{
{#5}_0\ar[d]\\
{#1}_0\ar[r]\ar@{=}[d] & \cdots \ar[r] & {#1}_{#4}\dequivarrow\\
{#1}\ar[rr] && {#3}\\ }
}

\def\labeldecomp#1#2#3#4#5#6#7#8{
\xymatrix{
{#5}_0 \ar[d]^{{#6}_0} && {#5}_1 \ar[d]^{{#6}_1} \\
{#1}_0\ar[rr]^{{#7}_0} \dequivarrow && {#1}_1\ar[rr]^{{#7}_1}&& \cdots 
\ar[rr]^{{#7}_{{#4}-1}} && {#1}_{#4}\ar@<-.4ex>[d]_{{#2}_{#4}}\\
X\ar[rrrrrr]^{#2} &&&&&& Y \ar@<-.4ex>[u]_{{#8}} \\ }
}

\def\shortlabeldecomp#1#2#3#4#5#6#7#8{
\xymatrix{
{#1}_0\ar[r]^{{#7}_0} \ar@{=}[d] & {#1}_1\ar[r]^{{#7}_1}& 
\cdots \ar[r]^(.42){{#7}_{{#4}-2}}
 & {#1}_{#4 -1 }\ar[r]^(.52){{#7}_{{#4}-1}} & {#1}_{#4}
\ar@<-.4ex>[d]_{{#2}_{#4}}\\
X\ar[rrrr]^{#2} &&&& Y \ar@<-.4ex>[u]_{{#8}} \\ }
}

\def\labeldecomp#1#2#3#4#5#6#7#8{
\xymatrix{
{#5}_0\ar[d]&{#5}_1\ar[d] && {#5}_{#4 -1}\ar[d]\\
{#1}_0\ar[r]^{{#7}_0} \ar@{=}[d] & {#1}_1\ar[r]^{{#7}_1}& 
\cdots \ar[r]^(.42){{#7}_{{#4}-2}}
 & {#1}_{#4 -1 }\ar[r]^(.52){{#7}_{{#4}-1}} & {#1}_{#4}
\ar@<-.4ex>[d]_{{#2}_{#4}}\\
X\ar[rrrr]^{#2} &&&& Y \ar@<-.4ex>[u]_{{#8}} \\ }
}

\begin{document}
\setlength{\abovedisplayskip}{\stdskip}     
\setlength{\belowdisplayskip}{\stdskip}     
%
%
%
%
%
%
%
%
\def\title#1{\def\thetitle{#1}}
\def\authors#1{\def\theauthors{#1}}
\def\author#1{\def\theauthors{#1}}
\def\address#1{\def\theaddress{#1}}
\def\secondaddress#1{\def\thesecondaddress{#1}}
\def\email#1{\def\theemail{#1}}
\def\url#1{\def\theurl{#1}}
\long\def\abstract#1\endabstract{\long\def\theabstract{#1}}
\def\primaryclass#1{\def\theprimaryclass{#1}}
\def\secondaryclass#1{\def\thesecondaryclass{#1}}
\def\keywords#1{\def\thekeywords{#1}}
%
\def\ifundefined#1{\expandafter\ifx\csname#1\endcsname\relax}
%
%
\long\def\maketitlepage{    

\vglue 0.2truein   

%
{\parskip=0pt\leftskip 0pt plus 1fil\def\\{\par\smallskip}{\LARGE
\bf\thetitle}\par\medskip}   

\vglue 0.15truein 

%
{\parskip=0pt\leftskip 0pt plus 1fil\def\\{\par}{\sc\theauthors}
\par\medskip}%
 
\vglue 0.1truein 

%
{\parskip=0pt\small
{\leftskip 0pt plus 1fil\def\\{\par}{\sl\theaddress}\par}
\ifundefined{thesecondaddress}\else\cl{and}  
{\leftskip 0pt plus 1fil\def\\{\par}{\sl\thesecondaddress}\par}\fi
\ifundefined{theemail}\else  
\vglue 5pt \def\\{\ \ {\rm and}\ \ } 
\cl{Email:\ \ \tt\theemail}\fi
\ifundefined{theurl}\else    
\vglue 5pt \def\\{\ \ {\rm and}\ \ } 
\cl{URL:\ \ \tt\theurl}\fi\par}

\vglue 7pt 

{\bf Abstract}

\vglue 5pt

\theabstract

\vglue 7pt 

{\bf AMS Classification numbers}\quad Primary:\quad \theprimaryclass

Secondary:\quad \thesecondaryclass

\vglue 5pt 

{\bf Keywords:}\quad \thekeywords


}    
%
%
%
%
%
%
%
%
%
%
%
%
%

\title{The Cone Length and Category of Maps:  Pushouts, Products and Fibrations}                    
\authors{Martin Arkowitz, Donald Stanley and Jeffrey Strom}                  
\address{}                  



\abstract   
For any collection of spaces $\A$, we investigate
two non-negative integer homotopy invariants of maps:  $\LL_\A(f)$,
the $\A$-cone length of $f$, and
$\L_\A(f)$, the $\A$-category of $f$.  When $\A$ is the collection of 
all spaces, these are the cone length and category of $f$, respectively,
both of which have been studied previously.  
The following results are obtained:
(1)  For a map of one homotopy pushout diagram into another, we derive
an upper bound for $\LL_\A$ and $\L_\A$ of the induced map of homotopy
pushouts in terms of $\LL_\A$ and $\L_\A$ of the other maps.  This has
many applications, including an inequality for $\LL_\A$ and $\L_\A$ of the
maps in a mapping of one mapping cone sequence into another.
(2)  We establish an upper bound for $\LL_\A$ and $\L_\A$ of the product
of two maps in terms of $\LL_\A$ and $\L_\A$ 
of the given maps and the $\A$-cone length of their domains.
(3) We study our invariants in a pullback square and obtain as a
consequence an upper bound for the $\A$-cone length and $\A$-category
of the total space of a fibration in terms of the 
$\A$-cone length and $\A$-category of the base and fiber.  
We conclude with several remarks, examples and open questions.
\endabstract


\primaryclass{55M30}                
\secondaryclass{55P99, 55R05}              
\keywords{Lusternik-Schnirelmann category, cone length, 
homotopy pushouts, products of maps, fibrations}                    

\maketitlepage
%
%

\begin{changemargin}{-1.1cm}{-1.1cm}

\section{Introduction}

In this paper we continue our investigation, begun in \cite{ASS}, of the cone length 
and category of   maps relative to a fixed collection of spaces.  For a  
collection $\A$ of spaces we consider 
two non-negative integer homotopy invariants of maps: the $\A$-category,
denoted  $\L_\A$,
and the $\A$-cone length, denoted $\LL_\A$.   When $\A$ is the collection of
all spaces, $\L_\A(f)$ is the category of the map $f$ as defined and studied in 
\cite{Fa-Hu}  and \cite{Cornea2} and $\LL_\A(f)$ is the
cone length of the map $f$ as defined and studied in \cite{Marcum} and \cite{Cornea1}.  
If, in this special case, $f$ is the inclusion of the base point
into $Y$, then $\L_\A(f)$ is just the category of the space $Y$,
which was introduced in 1934 by Lusternik and Schnirelmann in their work 
on the number of critical points of
smooth functions on a manifold \cite{L-S}.
In addition, $\LL_\A(f)$ is the cone length of $Y$ which 
has been studied by several people
\cite{Cornea1,Cornea4,Cornea5,Ganea,Marcum,Stanley,Takens} 
in the context of homotopy theory.  For an arbitrary collection
$\A$, the $\A$-category and the $\A$-cone length 
of the inclusion of the basepoint into $Y$ coincide  with
the $\A$-category and $\A$-cone length of $Y$.
Variants of this concept have been studied previously
\cite{S-T}. 

Thus our invariants are common generalizations of the category and cone length of a 
map and the $\A$-category and $\A$-cone length of a space.  In addition to 
providing a general framework for many
existing notions and retrieving known results as special cases, they have 
led to several new concepts and
results.  To discuss this, we first briefly summarize that part of our 
previous work which is relevant to this
paper.  More details are given in \sec 2.

In \cite{ASS} we introduced, for a fixed collection $\A$, five simple axioms which an 
integer valued function of based maps may satisfy.  
Then $\L_\A$ was defined as the maximum of
all such functions.  Similarly, $\LL_\A$ was defined as the maximum of all 
functions which satisfy an analogous
set of five axioms.  We then gave alternate characterizations of these 
invariants in terms of
certain decompositions of maps.  For instance, $\LL_\A(f)$ is essentially the smallest 
integer $n$ such that $f$ admits a decomposition up to homotopy as 
$$
X = X_0\maprt{j_0}X_1\maprt{j_1}\cdots \maprt{j_{n-1}} X_n = Y,
$$
where $L_i\maprt{}X_i\maprt{j_i}X_{i+1}$ is a mapping cone sequence with $L_i \in \A$.  
When $f$ is the inclusion of the base point into $Y$, we obtain $\cat_\A(Y)$ and $\cl_\A(Y)$, 
as noted above, and when $g$ is the map of $X$ to a  one point space, we obtain two new
invariants of spaces, the $\A$-kitegory of $X$, $\kit_\A(X) = \L_\A(g)$, and the 
$\A$-killing length of $X$, $\kl_\A(X) = \LL_\A(g)$.   In \cite{ASS} 
we made a preliminary study of these invariants and their interrelations.

From the time the concept of category was first introduced to the present, 
many people have been interested in the following questions:  
What is the relationship between the categories of the spaces 
which appear in a homotopy pushout \cite{Marcum,Hardie}?  
What is the category of the product of two spaces in terms of the categories 
of the factors \cite{A-Sta,Bassi,ClP2,Ganea2,Iwase1,Roitberg,Stanley2,Takens,Van}?  
What is the relationship of the categories of the spaces which appear
in a fiber sequence \cite{Hardie2,James-Singhof,V}?  
Similar questions have also been considered for cone length.  
In this paper we study these questions for the $\A$-category and $\A$-cone 
length of maps, and provide reasonably complete answers. 
This gives both new and known results for the $\A$-category and $\A$-cone length 
of spaces as well as new results for the $\A$-kitegory and the $\A$-killing length.

We now summarize the contents of the paper.  In \sec 2 we give our 
terminology and notation and discuss our earlier work in more detail.  
In \sec 3 we prove one of our main results, the Homotopy Pushout Mapping Theorem.  
This theorem gives an inequality for the $\A$-categories and
the $\A$-cone lengths of the four maps which constitute a map of one 
homotopy pushout square into another. 
Many applications  are given in \sec 4.  In particular, 
we obtain results about the
$\A$-categories and $\A$-cone lengths of the maps which appear in a mapping 
of one mapping cone sequence into
another.  In \sec 5 we establish an upper bound for the $\A$-category 
(resp., $\A$-cone length) of the product of
two maps in terms of the $\A$-categories (resp., $\A$-cone lengths) of the original 
maps and the $\A$-cone lengths of
their domains.  By specializing to spaces and letting $\A$ be the collection of 
all spaces, we retrieve classical
results on the category and cone length of the product of two spaces.  
We study pullbacks in \sec 6.  As a consequence of our main result on pullbacks
we obtain an inequality 
for the $\A$-category of the total space of a fiber sequence in terms of the
$\A$-category of the base and fiber, and a similar result for $\A$-cone length.  
Section 7 contains a potpouri of
results, examples and questions.   We begin by presenting a few simple, 
but useful, results about $\A$-category and
$\A$-cone length.  We then give some examples to illustrate the difference 
between these invariants for different
collections $\A$.  In particular, we show that some results that are known 
for the collection of all spaces  do
not hold for arbitrary collections.  Finally, we state and discuss some 
open problems.

We conclude this section by emphasizing two important points.
First of all, there are several different notions of the category of a map
in the literature.   The one that we generalize here to the $\A$-category 
of a map has been studied in \cite{Fa-Hu} and \cite{Cornea2}.  It is not the 
same as the one considered in \cite{Fox,B-G}.    In addition, Clapp and Puppe have 
considered the category of a map with respect to a collection of
spaces \cite{ClP,ClP2}.  However, their notion is completely different from ours.  
Secondly, although we state and prove our results in the category of well-pointed 
spaces and based maps, it should be clear that  nearly all our results hold in a 
(closed) model category \cite{Quillen} and that all of our results hold in a 
J-category \cite{Do,H-L}.

\bigskip

\section{Preliminaries}

In this section we give our notation and terminology and also recall some 
results from \cite{ASS} which will be needed later.  

All topological spaces are based and have the based homotopy type of 
CW-complexes, though we could more generally consider well-pointed based spaces.  
All maps and homotopies are to preserve base points.  We {\it do} distinguish 
between a map and a homotopy class.  By a commutative diagram we mean one which is
{\it strictly} commutative.  

We next give some notation which is standard for homotopy theory: 
$*$ denotes the base point of a space or the space consisting of a single point, 
$\simeq$ denotes homotopy of maps and $\equiv$ denotes same homotopy type of spaces. 
We let $0:X\maprt{} Y$ stand for the constant map and id$:X\maprt{} X$ for the 
identity map.   We use $\Sigma $ for (reduced) suspension, $*$ 
for (reduced) join, $\vee$  for wedge sum and $\smsh$ for smash product.  

We call a sequence $A\maprt{f} X\maprt{j} C$ of spaces and maps 
a  {\it mapping cone sequence} if
$C$  is the mapping cone of $f$ and $j$ is the standard inclusion.  
Then $j$ is a cofibration with cofiber $\Sigma A$.  
Using the mapping cylinder construction, we see that the  
concept of a cofiber sequence  and the concept of a mapping cone sequence are equivalent
\cite[Ch.\thinspace 3]{Hi}.  For maps $C\maplft{g}A\maprt{f}B$ we can form the homotopy pushout $Q$
$$
\xymatrix{
A\ar[rr]^f \ar[d]_g && B\ar[d] \\
C\ar[rr] && Q \\ }
$$
by defining $Q$ to be the quotient of $B\vee (A\times I) \vee C$ 
under the equivalence relation $(*,t)\sim *$, $(a,0) \sim f(a)$
and $(a,1) \sim g(a)$ for   $t\in I$ and $a \in A$.  
Note that $A \maprt{f} X \maprt{j} C$ is a mapping
cone sequence if and only if   
$$ 
\xymatrix{
A\ar[rr]^f \ar[d]  && X\ar[d]^j \\
{*}\ar[rr] && C \\ }
$$
is a homotopy pushout square.
The pullback $P$ of $C\maprt{g}A\maplft{f}B$ is defined by
$$
P = \{ (b,c)\, |\, b\in B, c\in C, f(b) = g(c)\} \subseteq B\times C. 
$$
We only use this construction when $f$ is a fibration.  
Thus all our pullbacks are homotopy pullbacks as well.  
Given a map $f:X\maprt{} Y$ we say that a map $g:X'\maprt{} Y'$ {\it homotopy dominates}
$f$ (or $f$ is a {\it homotopy retract} of $g$) if there is a homotopy-commutative 
diagram
$$
\xymatrix{
X\ar[rr]^i\ar[d]_f  &&   X'\ar[rr]^r\ar[d]^g   &&  X\ar[d]^f
\\
Y\ar[rr]^j        &&   Y'\ar[rr]^s         &&  Y 
}
$$
such that $ri \simeq \id$ and $sj \simeq \id$.  
If the diagram is strictly commutative and both homotopies
are equality, we delete the word `homotopy' from the definition.  
If $g$ homotopy dominates $f$ as above and in addition 
$ir\simeq \id$ and $js \simeq \id$ (i.e., $r$ and $s$ are homotopy equivalences 
with homotopy
inverses $i$ and $j$), we say that $f$ and $g$ are {\it homotopy equivalent}.

Next we recall some definitions and results from \cite{ASS} which will be used in the
sequel.   By a {\it collection} $\A$ we mean a class of spaces containing $*$ such that if
$A \in \A$ and $A\equiv A'$, then $ A' \in \A$.  We say that (1) $\A$ is {\it closed under
suspension} if 
$A \in \A $ implies $\Sigma A \in \A $, (2) $\A $ is {\it closed under wedges} 
if $A, A' \in \A $
implies $A\vee A' \in \A $ and (3) $\A $ is {\it closed under joins} if 
$A, A' \in \A $
implies $A *  A' \in \A $.  Examples of collections that we consider are 
(1) the collection $\A = \{{\mathrm{all\ spaces}}\}$ of all spaces, 
(2) the collection $\Sigma$ of all suspensions, and 
(3) the collection $\S$ of all wedges of spheres (including $S^0$).

Let $\A$ be a collection and $\l_{\A }$ an function which assigns to 
each map $f$  an integer
$0\le \l_{\A }(f) \le \infty $.  We say that $\l_{\A }$ satisfies the 
$\A$-{\it cone axioms} if 
\begin{enumerate}
\item[(1)] (Homotopy Axiom) If
$f \simeq g$, then $\l_{\A }(f)  = \l_{\A }(g)$.
\item[(2)]  (Normalization Axiom)  If
$f$ is a homotopy equivalence, then $\l_{\A }(f) =0$.
\item[(3)] 
(Composition Axiom)  $\l_{\A }(fg) \le  \l_{\A }(f) +  \l_{\A }(g)$.
\item[(4)]  
(Mapping Cone Axiom)  If $A\maprt{} X   \maprt{f} Y$ is a mapping cone sequence
with $A \in \A$, then $\l_{\A }(f) \le 1$. 
\item[(5)]  (Equivalence Axiom) 
If $f$ and $g$ are homotopy equivalent, then $\l_{\A }(f) = \l_{\A }(g) $.
\end{enumerate}
We say that $\l_{\A }$ satisfies the $\A $-{\it category axioms} if $\l_{\A}$ satifies
(1) -- (4) and 
\begin{enumerate}
\item[(5$'$)] (Domination Axiom) If $f$ is dominated by $g$, then $\l_{\A }(f) \le \l_{\A }(g) $.
\end{enumerate}

\bdfn{\rm
We denote by $\LL_{\A}(f)$ the maximum of all $\l_{\A }(f)$ where $\l_{\A }$ satisfies
(1)--(5) and by $\L_{\A}(f)$ the maximum of all $\l_{\A }(f)$ where $\l_{\A }$ satisfies
(1)--(4) and (5$'$).  We call $\LL_{\A}(f)$ the $\A $-{\it cone length} of $f$ and
$\L_{\A}(f)$ the $\A $-{\it category} of $f$.  }
\edfn

Since (5) is weaker than (5$'$), 
$\L_{\A}(f) \leq \LL_{\A}(f)$.   In  \cite{ASS} it is proved that when 
$\A = \{ {\mathrm{all\ spaces}} \}$, $\LL_{\A}(f)$ is
the cone length of $f$ as defined in \cite{Cornea2,Marcum}, 
and $\L_{\A}(f)$ is the category of $f$
as defined in \cite{Fa-Hu,Cornea2}.

One of the main results of \cite{ASS} gives alternate characterizations 
of $\LL_{\A}(f)$ and 
$\L_{\A}(f)$ in terms of decompositions of the map $f$.  
If $f:X\maprt{} Y$ is a map, then an {\it $\A $-cone decomposition of length $n$} 
of $f$  is
a homotopy-commutative diagram 
$$
\shortlabeldecomp{X}{f}{Y}{n}{K}{i}{j}{s}
$$
in which $f_n$ is a homotopy equivalence with homotopy inverse $s$  
and each map $j_i$ is part of a mapping cone sequence
$
\xymatrix@1{
A_i \ar[r] &  X_i \ar[r]^-{j_i}& X_{i+1}}
$
with $A_i\in \A$.  Thus $f_nj_{n-1}\cdots j_0
\simeq f$, $sf \simeq j_{n-1}\cdots j_0$,
$f_n s \simeq \id$ and $s f_n\simeq \id$. 
The 
homotopy-commutative diagram above is an
{\it $\A$-category  decomposition of $f$ of length $n$} 
if $s$ is simply a homotopy section of $f_n$, i.e., if  
$f_nj_{n-1}\cdots j_0 \simeq f$, $sf \simeq j_{n-1}\cdots j_0$ and
$f_n s \simeq \id$, but $sf_n$ need not be homotopic to the identity.
We prove in \cite[Thm.\thinspace 3.7]{ASS} that  
$$ \LL_{\A}(f) =\left\{ 
\begin{array}{ll}
0  & {\mathrm{if}}\ f\ {\mathrm {is\ a\ homotopy\ equivalence}}
\\
\infty & {\mathrm {if\ there\ is\ no}}\ \A{\mathrm{-cone\ decomposition\ of}}\ f
\\
n & {\mathrm {if}}\  n\ 
{\mathrm{is\ the\ smallest\ integer\ such\ that\ there\ exists\ an}}\\ & 
\A{\mathrm{-cone\ decomposition\ of\ length}}\  n\ {\mathrm{of}}\ f .\\
\end{array} \right.
$$
We similarly characterize $\L_{\A}(f)$ using $\A$-category decompositions instead of 
$\A$-cone decompositions.  Observe that if  the induced map $\pi_0(f):\pi_0(X) \to \pi_0(Y)$ on
path components is not surjective, then $L_\A(f)$ and $\L_\A(f)$ are
infinite for {\it every} collection $\A$. 

We have also studied four numerical invariants of 
spaces, defined in terms of the invariants $\LL_\A$ and $\L_\A$ as follows:
$$
\begin{array}{rclcrcl}  
\cl_{\A}(X) & = & \LL_\A(*\maprt{}X) &\quad &
\cat_{\A}(X) & = & \L_\A(*\maprt{}X)\\
      \kl_{\A}(X) & = & \LL_\A(X\maprt{}*) &\quad &
      \kit_{\A}(X) & = & \L_\A(X\maprt{}*) . \\ 
\end{array}
$$
When $\A =\{ {\mathrm{all\ spaces}}\}$, 
$\cat_\A(X) = \cat(X)$, the reduced Lusternik-Schnirelmann category of $X$
\cite[Prop.\thinspace 4.1]{ASS}, and 
$\cl_\A(X)= \cl(X)$, the cone length of $X$.
Moreover, $\kit_\A(X)\leq 1$ and $\kl_\A(X)\leq 1$ for every space $X$ in this case.

\bigskip

\section{The Homotopy Pushout Mapping Theorem}

In this section we prove the first main result of this paper.  
This consists of two inequalities, one for the $\A$-cone length 
and one for the $\A$-category of the maps from one homotopy pushout
square to another.   In Section \sec 4 
we will derive numerous consequences of this result.

We begin with a technical result that plays a key role in the proof of
the main theorem.

\blem\label{lem:strict}
Let $f: X\maprt{}Y$ be a  map with $\L_\A(f) = n$.  Then there exists
a map $g:X\maprt{} Z$ and maps $i:Y\maprt{}Z$ and $r:Z\maprt{}Y$
such that the diagram
$$
\xymatrix{   
           && X\ar[lld]_f\ar[rrd]^f \ar[d]_g         \\
Y\ar[rr]_i && Z\ar[rr]_r                        && Y,\\ 
}
$$
commutes, $ri= \id$ and $\LL_\A(g) = n$.
\elem

\begin{prf}
By \cite[Cor.\thinspace 4.4]{ASS} there is a homotopy-commutative diagram
$$
\xymatrix{
&&X \ar[lld]_f   \ar[d]_h \ar[rrd]^f\\
Y\ar[rr]_j && W \ar[rr]_s && Y \\
}
$$
with $sj\simeq\id$ and $\LL_\A(h) = n$.  We factor $s$ as
$W\maprt{s^0} E\maprt{s^1}Y$, where $s^0$ is a homotopy equivalence
and $s^1$ is fibration \cite[Ch.\thinspace 3]{Hi}.  Then $s^0j :Y\maprt{}E$ and $s^1s^0j\simeq \id$.
Thus there is a map $t:Y\maprt{}E$ such that 
$t\simeq s^0 j$ and $s^1t =\id$.  Now we have a homotopy-commutative diagram
$$
\xymatrix{
&&X  \ar[lld]_f  \ar[d]_{s^0h} \ar[rrd]^f   \\
Y\ar[rr]_t            && E \ar[rr]_{s^1}             && Y.          \\
}
$$
We factor $s^0$ as $W\maprt{k'}E'\maprt{k''}E$
where $k'$ is a homotopy equivalence and $k''$ is a fibration.
Then $k'h:X\maprt{}E'$ and $k''k'h =s^0h\simeq tf$, and so
there is a map $l:X\maprt{}E'$ with $l\simeq k'h$ and $k''l = tf$.
Then we have a diagram
$$
\xymatrix{
&& X  \ar[lld]_f   \ar[d]_{k''l} \ar[rrd]^f\\
Y\ar[rr]_t && E \ar[rr]_{s^1} && Y \\
}
$$
with $s^1t =\id$ and the left hand triangle strictly commutative.
It then follows that the right hand triangle is strictly commutative.
Note that $k''l\simeq k''k'h =s^0 h$.  Since $s^0$ is a homotopy
equivalence, $\LL_\A(k''l) = \LL_\A(h)$.  This proves the lemma
with $Z= E$, $g= k''l$, $i=t$ and $r=s^1$.
\end{prf}

\bthm\label{thrm:MPT}
Let $\A$ be a collection of spaces that is closed under wedges and
suspension and let
$$
\xymatrix{
C\ar[d]^{c} && A \ar[ll]_{g} \ar[rr]^{f} \ar[d]^{a} && B\ar[d]^{b}   
\\
C' && A' \ar[ll]_{g'}\ar[rr]^{f'} && B'   
\\  }
$$
be a commutative diagram.  Let $D$ be the homotopy pushout of 
the top row, $D'$ be the homotopy pushout of the bottom row,
and $d:D\maprt{}D'$ the induced map.  Then
\begin{enumerate}
\item  $\LL_\A(d) \leq \LL_\A(a) + 
\max ( \LL_\A(b),\LL_\A(c))$;
\item  $\L_\A(d) \leq \L_\A(a) + 
\max ( \L_\A(b),\L_\A(c))$.
\end{enumerate}
\ethm

\brmk{\rm
Our proof will show that
if $\A$ is only assumed to be closed under wedges, then (1) and (2) hold when 
$A=A'$ and $a=\id$.  
In Corollary \ref{cor:nowedge} below we derive some slightly weaker formulas
that require only closure under suspension. 
In particular, if $\A$ is only known to be closed under 
suspensions, then (1) and (2) hold with $B= B'$, $C= C'$, $b=\id$ and $c=\id$.}
\ermk

\begin{prf}
First we prove (1).  We factor the given diagram as
$$
\xymatrix{
C\ar[d]_c && A \ar[ll]_g\ar[rr]^f\ar@{=}[d]  && B\ar[d]^b   \\
C' \ar@{=}[d] && A \ar[ll]_{cg}\ar[rr]^{bf}\ar[d]^{a} && B'\ar@{=}[d]   \\ 
C' && A' \ar[ll]_{g'}\ar[rr]^{f'} && B'    \\  }
$$
and let $\overline D$ denote the homotopy pushout of the middle row.
Then we have a factorization 
$$
\xymatrix{ 
D\ar[rd]_{d'} \ar[rr]^d && D' \\
&\overline D \ar[ru]_{d''}\\ }
$$
of $d$.
Since $\LL_\A(d) \leq \LL_\A(d') + \LL_\A(d'')$ by the Composition Axiom, 
it suffices to prove the result in the two special cases 
\begin{enumerate}
\item[(a)]  $A=A'$ and $a=\id$,
\item[(b)] $B=B'$, $C=C'$, $b=\id$ and $c=\id$.
\end{enumerate}

We begin with (a) and let $m = \max ( \LL_\A(b),\LL_\A(c))$.
We consider an $\A$-cone decomposition of $b$ of length
$m$.  This yields a homotopy factorization of 
$b\simeq hi_{m-1}\cdots i_1i_0$:
$$
B = X_0\maprt{i_0}X_1\maprt{i_1}\cdots \maprt{i_{m-1}}X_m\maprt{h}B',
$$
where $A_l\maprt{}X_l\maprt{i_l}X_{l+1}$ is a mapping cone sequence 
with $A_l\in \A$ for each $l$ and $h$ is a homotopy equivalence.
Since $i_{m-1}\cdots i_1i_0$ is a cofibration, $h$ is homotopic
to a map (also called $h$) such that $b = hi_{m-1}\cdots i_1i_0$.
Similarly, we have an $\A$-cone decomposition of $c$ of length
$m$ which gives a factorization $c = k j_{m-1}\cdots j_1 j_0$:
$$
C = Y_0\maprt{j_0}Y_1\maprt{j_1}\cdots \maprt{j_{m-1}}Y_m\maprt{k}C',
$$
where $B_l\maprt{} Y_l\maprt{j_{l}}Y_{l+1}$ is a mapping cone sequence
with $B_l\in\A$ for each $l$ and $k$ is a homotopy equivalence.  Thus we 
have a commutative diagram
$$
\xymatrix{
Y_0=C \ar[d]_-{j_0} && A \ar[ll]_-g\ar[rr]^-{f}\ar@{=}[d] && B=X_0 \ar[d]^-{i_0}   
\\
Y_1\ar[d]_-{j_1}   && A \ar[ll]_-{j_0g}\ar[rr]^-{i_0f}\ar@{=}[d] && X_1 \ar[d]^-{i_1}   
\\ 
\vdots \ar[d] && \vdots \ar@{=}[d] && \vdots \ar[d] 
\\
Y_{m-1} \ar[d]_-{j_{m-1}} && A 
\ar[ll]_-{j_{m-2} \cdots j_0 g}
\ar[rr]^-{i_{m-2} \cdots i_0 f}
\ar@{=}[d] && X_{m-1} \ar[d]^-{i_{m-1}}   
\\
Y_{m} \ar[d]_k   && A 
\ar[ll]_-{j_{m-1} \cdots j_0 g}
\ar[rr]^{i_{m-1} \cdots i_0 f}
\ar@{=}[d] && X_{m} \ar[d]^-{h}    
\\ 
C' && A  \ar[ll]_-{cg}\ar[rr]^-{bf} &&   B' .  
\\  }
$$
We number the rows $0, 1, \ldots, m+1$ and let $D_l$ be the homotopy
pushout of the $l^{th}$ row, with induced maps $d_l:D_l\maprt{} D_{l+1}$.
Then $D_0 = D$, $D_{m+1} = \overline D$ and $d_m\cdots d_0 = d':D\maprt{} \overline D$.
Thus it suffices to prove 
\begin{enumerate}
\item[(i)]  $\LL_\A(d_l) \leq 1$ for $l= 0, \ldots, m-1$,
\item[(ii)]  $\LL_\A(d_m) =0$.
\end{enumerate}
We first establish (i).  Consider the commutative diagram
$$
\xymatrix{
B_l\ar[d] && {*}\ar[ll]\ar[rr]\ar[d]  && A_l\ar[d]  &&&& A_l\wdg B_l\ar[d] \\
Y_l\ar[d]_{j_l} && A \ar[ll]\ar[rr]\ar@{=}[d] && X_l\ar[d]^{i_l} &
{}\ar@{.>}[rr]^{ {\mathrm{homotopy}}\atop {\mathrm{pushout}} } &&&  
  D_l\ar[d] \\
Y_{l+1}      && A\ar[ll]\ar[rr]         && X_{l+1}  &&&& D_{l+1},      \\   } 
$$
where the columns are regarded as mapping cone sequences.
The homotopy pushouts of the rows form a sequence
$A_l\wdg B_l \maprt{} D_l\maprt{} D_{l+1}$.  By the Four
Cofibrations Theorem, this is a cofiber sequence (see \cite[p.\thinspace 21]{Do2}).
Since $\A$ is closed under wedges, $A_l\wdg B_l\in\A$.  
Therefore $\LL_\A(d_l)\leq 1$ by the Mapping Cone Axiom.  For (ii) we note that
$d_m: D_m\maprt{}\overline D$ is a homotopy equivalence since
$h$ and $k$ are homotopy equivalences \cite[Ch.\thinspace XII,\thinspace \sec\thinspace 4.2]{B-K}.  Thus
$\LL_\A(d_m) =0$, which completes the proof of (a).

For (b) we proceed similarly by assuming that $\LL_\A(a) = m$ and taking
an $\A$-cone length decomposition of $a$ of length $m$:
$$
A = X_0\maprt{i_0} X_1\maprt{} \cdots\maprt{} X_{m-1}\maprt{i_{m-1}} X_m \maprt{h} A',
$$
where $A_l \maprt{} X_l\maprt{i_l} X_{l+1}$ is a mapping cone sequence
with $A_l\in \A$, $h$ a homotopy equivalence and $a= hi_{m-1}\cdots i_0$.
This yields a commutative diagram
$$
\xymatrix{
C'                                         \ar@{=}[d] 
&& 
A=X_0 
\ar[ll]_-{cg}\ar[rr]^-{bf} \ar[d]^-{i_{0}} 
&& B'                                      \ar@{=}[d]   
\\
C'                                         \ar@{=}[d]   
&& 
X_1
\ar[ll]_-{g'hi_{m-1}\cdots i_1 }
\ar[rr]^-{f'hi_{m-1}\cdots i_1 }
\ar[d]^-{i_{1}} 
&& 
B' \ar@{=}[d]
\\ 
\vdots \ar@{=}[d]  && \vdots  \ar[d]  && \vdots    \ar@{=}[d] \\
C'                                         \ar@{=}[d] 
&& 
X_{m-1}  
\ar[ll]_-{g'hi_{m-1}} 
\ar[rr]^-{f'hi_{m-1}} 
\ar[d]^-{i_{m-1}} 
&& 
B'                                          \ar@{=}[d]   
\\ 
C'                                          \ar@{=}[d]  
&&X_m 
\ar[ll]_-{g'h}
\ar[rr]^-{f'h}
\ar[d]^-{h} 
&& B'                                       \ar@{=}[d]    \\ 
C' 
&& 
A'  
\ar[ll]_-{g'}\ar[rr]^-{f'} 
&&   B' .  \\  }
$$
We number the rows $0,1,\ldots, m+1$ and 
let $\twdl D_l$ be the homotopy pushout of the $l^{th}$ row
with induced maps $\twdl d_l:\twdl D_l\maprt{} \twdl D_{l+1}$.  Then 
$\twdl D_0 = \overline D$, $\twdl D_{m+1}= D'$ and 
$d'' = \twdl d_m\cdots \twdl d_1\twdl d_0$.  It suffices to show (i)
$\LL_\A(\twdl d_l) \leq 1$ for $l= 0, \ldots , m-1$ and (ii) $\LL_\A(\twdl d_m) =0$.
The argument is similar to (a), and so we content ourselves with noting 
that 
$$
\xymatrix{
{*}\ar[d]    &&   A_l\ar[ll]\ar[rr]\ar[d]       &&   {*}\ar[d]          &&&& \s A_l\ar[d] 
\\
C'\ar@{=}[d] &&   X_l\ar[ll]\ar[rr]\ar[d]^{i_l} && B'\ar@{=}[d] 
& {}\ar@{.>}[rr]^{  {\mathrm{homotopy}}\atop {\mathrm{pushout}} }        &&& \twdl D_l\ar[d]^{\twdl d_l} 
\\ 
C'           && X_{l+1} \ar[ll]\ar[rr]          && B'                   &&&& \twdl D_{l+1}     
\\  }
$$
determines a sequence $\s A_l \maprt{} \twdl D_l \maprt{\twdl d_l}D_{l+1}$
since $\s A_l$ is the homotopy pushout of the top row.  By the 
Four Cofibrations Theorem this is a cofiber sequence.  Since $\A$
is closed under suspension, $\s A_l\in\A$, and so $\LL_\A(\twdl d_l)\leq 1$.  
This completes the proof of (1).

To prove (2), we apply Lemma \ref{lem:strict}.  Thus, there are 
commutative diagrams
$$
\xymatrix{
& A\ar[rd]^a\ar[ld]_a \ar[d]_\alpha  \\
A'\ar[r]^i & X \ar[r] ^r & A' \\
}\quad
\xymatrix{
& B\ar[rd]^b\ar[ld]_b \ar[d]_\beta  \\
B'\ar[r]^j & Y \ar[r] ^s & B' \\
}\quad \xymatrix{{}
\ar@{}[d]^{\txt{and}}\\{}}\quad 
\xymatrix{
& C\ar[rd]^c\ar[ld]_c \ar[d]_\gamma  \\
C'\ar[r]^k & Z \ar[r] ^t & C' \\
}
$$
with $ri=\id$, $sj=\id$, and $tk =\id$, and $\LL_\A(\alpha) = \L_\A(a)$,
$\LL_\A(\beta) = \L_\A(b)$, and $\LL_\A(\gamma) =\L_\A(c)$.  Thus we have
a diagram
$$
\xymatrix{
& C \ar[rdd]^(.35)c  \ar[dd]^(.35)\gamma  \ar[ldd]_(.35)c 
\\
&&&& A \ar[lllu]_g\ar[rrrd]^f\ar[rdd]^(.35)a\ar[ldd]_(.35)a\ar[dd]^(.35)\alpha
\\
C'\ar[r]^k &Z\ar[r]^t & C'  &&&&& B\ar[rdd]^(.35)b\ar[ldd]_(.35)b\ar[dd]^(.35)\beta
\\
&&&A'\ar[r]^(.3)i \ar[lllu]^(.65){g'}\ar[rrrd]_(.3){f'} & X
\ar[r]^(.35)r\ar[lllu]_(.6)v\ar[rrrd]^(.4)u & A'\ar[lllu]_(.65){g'}\ar[rrrd]^(.35){f'}
\\
&&&&&&B'\ar[r]_j & Y\ar[r]_s & B'
\\
}
$$
where $u=jf'r:X\maprt{}Y$ and $v=kg'r:X\maprt{}Z$.
All triangles and rectangles in the above diagram are commutative.
If we denote by $E$ the homotopy pushout of $Z\maplft{v} X\maprt{u}Y$
and the induced maps of homotopy pushouts by $e:D\maprt{}E$, 
$l:D'\maprt{} E$ and $m:E\maprt{}D'$, then we have a commutative 
diagram
$$
\xymatrix{
&& D\ar[rrd]^d\ar[lld]_d \ar[d]_e \\
D'\ar[rr]^l && E \ar[rr] ^m && D' \\
}
$$
with $ml = \id$.  Therefore
$$
\begin{array}{rclcl}
\L_\A(d)  & \leq  & \L_\A(e)    &\quad &  {\mathrm{since}}\ e\ {\mathrm{dominates}}\ d \\
          & \leq  & \LL_\A(e) \\
          & \leq & \LL_\A(\alpha) + \max(L_\A(\beta),L_\A(\gamma) ) && {\mathrm{by\ part}}\ (1) \\
          & =    & \L_\A(a) + \max(\L_\A(b),\L_\A(c)).\\
\end{array}
$$
\end{prf}

We next show that
some of Theorem \ref{thrm:MPT} remains true with weaker hypotheses on the 
collection $\A$.

\bcor\label{cor:nowedge}
Assume the hypotheses of Theorem \ref{thrm:MPT},
except that $\A$ is not necessarily closed under wedges.
Then
\begin{enumerate}
\item  $\LL_\A(d) \leq \LL_\A(a) + \LL_\A(b) + \LL_\A(c)$;
\item  $\L_\A(d) \leq \L_\A(a) + \L_\A(b) + \L_\A(c)$.
\end{enumerate}
The result remains true without assuming that 
$\A$ is closed under suspensions if $A = A'$ and
$a =\id$.
\ecor

\begin{prf}
We simply decompose the given map of homotopy pushouts into 
a composition of three maps:
$$
\xymatrix{
C  \ar[d]_c   && A \ar[ll]_g\ar[rr]^f\ar@{=}[d]       && B\ar@{=}[d]         \\
C' \ar@{=}[d] && A \ar[ll]_{cg}\ar[rr]^f\ar@{=}[d]       && B\ar[d]^b           \\
C' \ar@{=}[d] && A \ar[ll]_{cg}\ar[rr]^{bf}\ar[d]^{a} && B'\ar@{=}[d]        \\ 
C'            && A' \ar[ll]_{g'}\ar[rr]^{f'}          && B'  .                \\  }
$$
The method of the proof of Theorem \ref{thrm:MPT} is then applied to each
factor.
\end{prf}

\section{Applications of the Homotopy Pushout Mapping Theorem}

In this section we illustrate the power of the homotopy pushout mapping theorem by
obtaining as a consequence a large number of results, some known (in the case 
$\A = \{ {\mathrm {all\ spaces}} \}$), and some new.

\subsection{Homotopy Pushouts}

\bcor\label{cor:pushouts}  
Let $\A$ be any collection of spaces.  Let 
$$
\xymatrix{ A\ar[rr] \ar[d] && B\ar[d]\\ C\ar[rr] && D\\ }
$$
be a homotopy pushout square.  Then  
\begin{enumerate}

\item 
\begin{enumerate}
\item  $\LL_\A(B\maprt{}D) \leq \LL_\A(A\maprt{}C)$,
\item  $\L_\A(B\maprt{}D) \leq \L_\A(A\maprt{}C)$;  
\end{enumerate}

\item  
\begin{enumerate}
\item  $\cl_\A( D ) \leq \cl_\A(B) + \LL_\A(A\maprt{}C)$,
\item  $\cat_\A( D ) \leq \cat_\A(B) + \L_\A(A\maprt{}C)$;  
\end{enumerate}

\item
\begin{enumerate}
\item $\kl_\A(B) \leq \LL_\A(A\maprt{}C) + \kl_\A(D)$,
\item $\kit_\A(B) \leq \L_\A(A\maprt{}C) + \kit_\A(D)$.
\end{enumerate}

\item  If $A$ is closed under wedges, then
\begin{enumerate}
\item  $\LL_\A( A\maprt{}D)\leq \max( \LL_\A(A\maprt{}B) , \LL_\A(A\maprt{}C) )$,
\item  $\L_\A( A\maprt{}D)\leq \max( \L_\A(A\maprt{}B) , \L_\A(A\maprt{}C) )$.
\end{enumerate}

\end{enumerate}
\ecor

\begin{prf}
The proof of each part amounts to constructing the 
correct diagram.

\noindent{\sc Proof of 1}\quad   Apply Corollary \ref{cor:nowedge}   
to the diagram
$$
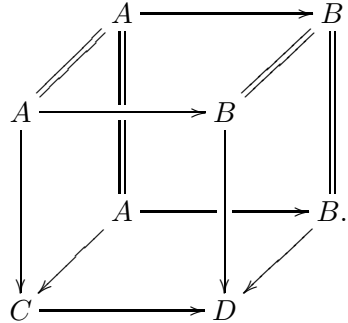
\xymatrix{
& A \ar[rr]\ar@{=}[ld]\ar@{=}'[d][dd]      && B \ar@{=}[dd] \ar@{=}[ld]  
\\
A \ar[rr] \ar[dd]                  && B \ar[dd] 
\\
& A\ar'[r][rr]\ar[ld]             && B. \ar[ld]
\\
C\ar[rr]                          && D
}
$$

\noindent{\sc Proof of 2 and 3} \quad Apply (1) to the 
diagram
$$
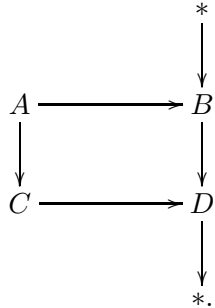
\xymatrix{ && {*}\ar[d] \\A\ar[rr] \ar[d] && B\ar[d]\\ C\ar[rr] && D\ar[d] \\
&& {*}.\\  }
$$

\noindent{\sc Proof of 4}\quad 
Map the trivial homotopy pushout diagram
$$
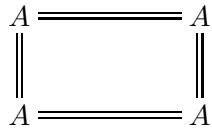
\xymatrix{ A\ar@{=}[rr] \ar@{=}[d] &&A\ar@{=}[d]\\ A\ar@{=}[rr] && A\\ }
$$
into the given one, and apply Theorem \ref{thrm:MPT}.
\end{prf}

\bcor\label{cor:susppushouts}
Let $\A$ be a collection of spaces that is closed under 
wedges and suspension and let
$$
\xymatrix{ A\ar[rr] \ar[d] && B\ar[d]\\ C\ar[rr] && D\\ }
$$
be a homotopy pushout square.  Then
\begin{enumerate}
\item
\begin{enumerate}
\item  $\cl_\A(D) \leq \cl_\A(A) + \max (\cl_\A(B), \cl_\A(C) )$,
\item  $\cat_\A(D) \leq \cat_\A(A) + \max (\cat_\A(B), \cat_\A(C) )$;
\end{enumerate}
\item
\begin{enumerate}
\item  $\kl_\A(D) \leq \kl_\A(A) + \max (\kl_\A(B), \kl_\A(C) )$,
\item  $\kit_\A(D) \leq \kit_\A(A) + \max (\kit_\A(B), \kit_\A(C) )$.
\end{enumerate}
\end{enumerate}
\ecor

\begin{prf}
For (1), apply Theorem \ref{thrm:MPT} to the map of the trivial homotopy
pushout diagram
$$
\xymatrix{ {*}\ar@{=}[rr] \ar@{=}[d] &&{*}\ar@{=}[d]\\ {*}\ar@{=}[rr] && {*}\\ }
$$
into the given homotopy pushout; for (2), map the given homotopy pushout into the
trivial one.
\end{prf}

\brmk
{\rm In the special case $\A=\{ {\mathrm{all\ spaces}} \}$,
Marcum \cite{Marcum} has proved  Corollary \ref{cor:pushouts}(1a)
and Hardie \cite{Hardie} has proved Corollary \ref{cor:susppushouts}(1b)
(see also \cite{Cornea4}).}
\ermk

\medskip

\subsection{Mapping Cone Sequences}

As noted in \sec 2, a mapping cone sequence $A\maprt{}B\maprt{}C$ 
can be regarded as a homotopy pushout square.  
Therefore the results of 4.1 apply to mapping cone sequences.

\bcor\label{cor:mappingcones}
Let $\A$ be any collection of spaces.  
Let $A\maprt{}B\maprt{}C$ be a mapping cone sequence.
Then
\begin{enumerate}

\item  
\begin{enumerate}
\item  $\cl_\A(C) \leq \LL_\A(A\maprt{} B )$,
\item  $\cat_\A(C ) \leq \L_\A(A\maprt{} B)$;
\end{enumerate}

\item
\begin{enumerate}
\item   $\LL_\A (B\maprt{}  C) \leq \kl_\A(A)$,
\item   $\L_\A(B\maprt{}  C) \leq \kit_\A(A)$;
\end{enumerate}
\item
\begin{enumerate}
\item  $\cl_\A(C) \leq \kl_\A(A) + \cl_\A(B)$,
\item  $\cat_\A(C) \leq \kit_\A(A) + \cat_\A(B)$;
\end{enumerate}

\item    
\begin{enumerate}
\item  $\kl_\A(B)\leq \kl_\A(A) + \kl_\A(C)$,
\item  $\kit_\A(B)\leq \kit_\A(A) + \kit_\A(C)$.
\end{enumerate}

\end{enumerate}
\ecor

\begin{prf}
\noindent{\sc Proof of 1 and 2}\quad  Immediate from Corollary \ref{cor:pushouts}(1).

\noindent{\sc Proof of 3 and 4} \quad Immediate from (2) and (3)
of Corollary \ref{cor:pushouts}.
\end{prf}

\brmk{\rm
Corollary \ref{cor:mappingcones}(4)   
shows that $\kl_\A$ and $\kit_\A$ are subadditive on cofibrations in the
following sense (we only state this for $\kl_\A$): If 
$A\maprt{}X\maprt{}Q$ is a cofiber sequence, then $\kl_\A(X) 
\leq \kl_\A(A) + \kl_\A(Q)$ (see \cite[Thm.\thinspace 3.4]{AS}.  
This follows (when $\A$ is closed under wedges) since
every cofiber sequence is equivalent to a mapping cone sequence.  This inequality is not
generally true for $\cl_\A$ or $\cat_\A$ as the cofiber sequence 
$$
S^2\maprt{}\cp^3\maprt{}S^4 \wdg S^6
$$
shows for the collections $\A = \S$, $\Sigma$ and $\{ {\mathrm{all\ spaces}} \}$.
}
\ermk

\bcor\label{cor:suspcof} Let $\A$ be a collection of 
spaces that is closed under suspension.
Consider the map of one mapping cone sequence into another given
by the commutative diagram  
$$
\xymatrix{
A \ar[rr]\ar[d] && B \ar[rr]\ar[d] && C\ar[d]\\
A'\ar[rr] && B'\ar[rr] && C'. \\ }
$$
Then
\begin{enumerate}
\item  $\LL_\A(C\maprt{} C') \leq \LL_\A(A\maprt{} A') + \LL_\A(B\maprt{}B')$,
\item  $\L_\A(C\maprt{} C') \leq \L_\A(A\maprt{} A') + \L_\A(B\maprt{}B')$.
\end{enumerate}
\ecor

\begin{prf}
Apply Corollary \ref{cor:nowedge}
to the homotopy pushouts obtained from the mapping cone sequences.
\end{prf}

\medskip

\subsection{Other Consequences}

\bcor\label{cor:klcl}
Let $\A$ be any collection of spaces.  Then
for any space $B$,
\begin{enumerate}
\item  $\cl_\A(\s B) \leq \kl_\A(B)$;
\item  $\cat_\A(\s B) \leq \kit_\A(B)$.
\end{enumerate}
\ecor

\begin{prf}  Apply  Corollary 
\ref{cor:pushouts}(1) to the homotopy pushout
square
$$
\xymatrix{B\ar[rr] \ar[d] &&{*}\ar[d]\\ {*}\ar[rr] && \s B.\\ }
$$
\end{prf}

\bcor\label{cor:4.10} 
Let $\A$ be a  collection of spaces that is closed under suspension.  
\begin{enumerate}
\item  For any map $f:A\maprt{}B$,  
\begin{enumerate}
\item  $\LL_\A(f) \leq \cl_\A(A) + \cl_\A(B) $,
\item  $\L_\A(f) \leq \cat_\A(A) + \cat_\A(B)$;
\end{enumerate}

\item  For any space $A$,  
\begin{enumerate}
\item  $\kl_\A( A) \leq \cl_\A(A)$,
\item  $\kit_\A( A) \leq \cat_\A(A)$;
\end{enumerate}

\item  If $f:A\maprt{} B$ and $g:B\maprt{}C$, then 
\begin{enumerate}
\item  $\LL_\A(g) \leq \LL_\A(f) + \LL_\A(g  f) $,
\item  $\L_\A(g) \leq \L_\A(f) + \L_\A(g  f) $;
\end{enumerate}

\item  If $f:A\maprt{}B$ and $g: B\maprt{}A$ with $gf =\id$, then 
\begin{enumerate}
\item[] $\L_\A(g) \leq \cat_\A(B)$;
\end{enumerate}

\item  If $f:A\maprt{} B$ and $g:B\maprt{}A$ with 
$g  f = \id$, then 
\begin{enumerate}
\item  $\LL_\A(g) \leq \LL_\A(f)  $,
\item  $\L_\A(g) \leq \L_\A(f) $.
\end{enumerate}

\end{enumerate}
\ecor

\begin{prf}
Again, the proofs depend on finding the appropriate diagram.

\noindent{\sc Proof of 1}\quad   Apply Corollary \ref{cor:nowedge}  
to the diagram
$$
\xymatrix{
& {*} \ar@{=}[rr]\ar[ld]\ar'[d][dd]      && {*} \ar[dd] \ar[ld]  
\\
A \ar@{=}[rr] \ar@{=}[dd]                  && A \ar[dd] 
\\
& A\ar'[r][rr]\ar@{=}[ld]             && B. \ar@{=}[ld]
\\
A\ar[rr]                          && B 
}
$$

\noindent{\sc Proof of 2}\quad   Apply (1) to the map $A\maprt{} *$.

\noindent{\sc Proof of 3}\quad   Apply Corollary \ref{cor:nowedge}
to the diagram
$$
\xymatrix{
& A \ar[rr]\ar@{=}[ld]\ar'[d][dd]      && B \ar@{=}[dd] \ar@{=}[ld]  
\\
A \ar[rr] \ar[dd]                  && B \ar[dd] 
\\
& B\ar@{=}'[r][rr]\ar[ld]             && B. \ar[ld]
\\
C\ar@{=}[rr]                          && C 
}
$$

\noindent{\sc Proof of 4}\quad   We consider the following mapping of
homotopy pushout squares  
$$
\xymatrix{
& {*} \ar[rr]\ar[ld] \ar'[d][dd]      && B \ar@{=}[dd] \ar[ld]  
\\
A \ar[rr] \ar@{=}[dd]                   && A\wdg B \ar[dd]^(.3){(\id,g)}
\\
& B\ar@{=}'[r][rr]\ar[ld]^(.4)g             && B. \ar[ld]^(.4)g
\\
A\ar@{=}[rr]                          && A 
}
$$
By Corollary \ref{cor:nowedge}, we immediately conclude that 
$\L(\id,g) \leq \cat (B)$.  Now the commutative diagram
$$
\xymatrix{
B\ar[d]_g  \ar[rr]    && A\wdg B \ar[d]^{(\id,g)} \ar[rr]^-{(f,\id)}&& B \ar[d]^g
\\
A\ar@{=}[rr]                   && A\ar@{=}[rr]                     && A  
\\
}
$$
shows that $g$ is dominated by $(\id, g)$.  Thus 
$\L(g) \leq \L(\id ,g) \leq \cat(B)$.

\noindent{\sc Proof of 5}\quad   Apply (3), 
using the fact that $L(\id) = \L(\id ) =0$.

\end{prf}

\bcor\label{cor:lowerbound}
Let $\A$ be a collection that is closed under  
suspension and let $f:A\maprt{}B$.  Then
\begin{enumerate}
\item
\begin{enumerate}
\item   $\LL_\A(f) \geq | \kl_\A(B) - \kl_\A(A) |$, 
\item   $\L_\A(f) \geq | \kit_\A(B) - \kit_\A(A) |$;
\end{enumerate}
\item
\begin{enumerate}
\item   $\LL_\A(f) \geq   \cl_\A(B) - \cl_\A(A)  $, 
\item   $\L_\A(f) \geq   \cat_\A(B) - \cat_\A(A)  $.
\end{enumerate}
\end{enumerate}
\ecor

\begin{prf}
We only prove (1a); the other parts are similar.  
The Composition Axiom, applied to $A\maprt{f} B\maprt{} *$, implies that
$$
\kl_\A(A) =    \LL_\A(A\maprt{}*) 
          \leq \LL_\A(B\maprt{}*) + \LL_\A(f) 
          =    \kl_\A(B) + \LL_\A(f),
$$
so $\LL_\A(f) \geq \kl_\A(A) - \kl_\A(B)$.  On the other hand, 
Corollary \ref{cor:4.10}(3) shows that
$$
\kl_\A(B) = \LL_\A(B\maprt{}*) \leq  \LL_\A(A\maprt{}*) + \LL_\A(f)  = \kl_\A(A) + \LL_\A(f),
$$
so $\LL_\A(f) \geq \kl_\A(B) - \kl_\A(A)$. This proves (1a)
\end{prf}

\brmk  {\em  Assume $\A$ is closed under suspension.  Then
by Corollary \ref{cor:4.10}, $\kl_\A(A) \leq \cl_\A(A)$ and $\kit_\A(A)\leq \cat_\A(A)$ 
(the first inequality was also proved in \cite[Thm.\thinspace 3.3]{AS}).
Corollary \ref{cor:klcl} then shows that $\cl_\A$ and $\kl_\A$
agree stably (this was proved by Christensen in \cite{Christensen}), 
and similarly for $\cat_\A$ and $\kit_\A$.  Additionally, Cornea \cite{Cornea2} 
has given a completely different proof of Corollary \ref{cor:4.10}(4) in the case
$\A = \{ {\mathrm{all\ spaces}} \}$. }
\ermk

\medskip

\subsection{Partial Converse to Theorem \ref{thrm:MPT}}

In this section we show that the formulas of 
Theorem \ref{thrm:MPT} very nearly characterize
those collections $\A$ which are closed under 
wedges or under suspensions.

We introduce the following new 
construction:  for any collection $\A$,  The collection 
$\overline \A$ is defined to be 
$$
\overline \A = \{ X \, |\, \kl_\A(X) \leq 1\}.
$$
Our first result shows that passing from $\A$ to $\overline \A$
has no effect on the corresponding cone length and 
category invariants.  Note that if every space in 
$\A$ is simply-connected, then
$\A = \overline \A $.

\bprp
For any map $f:X\maprt{}Y$, 
\begin{enumerate}
\item  $\LL_{\overline \A} (f) = \LL_\A(f)$,
\item  $\L_{\overline \A} (f) = \L_\A(f)$.
\end{enumerate}
\eprp

\begin{prf}
If suffices to prove (1), because for any collection $\A$,
$\L_\A(f)$ is the least $n$ for which $f$ is a retract of a map 
$g$ with $\LL_\A(g) \leq n$ \cite[Prop.\thinspace 4.3]{ASS}.  

Since $\A\sseq \overline \A$, we have $\LL_{\overline \A} (f) \leq \LL_\A(f)$
for any map $f$, so it remains to prove the reverse inequality.  
Suppose $\LL_{\overline \A} (f) = n$, and that 
$$
\xymatrix{ 
A_0  \ar[d]           & A_1\ar[d]  &                & A_{n-1}\ar[d]
\\
X_0 \ar[r]^-{j_0}\ar@{=}[d]  & X_1\ar[r]^-{j_1}  & \cdots\ar[r]^-{j_{n-2}}   & X_{n-1} \ar[r]^-{j_{n-1}}  &
X_n\dequivarrow 
\\
X\ar[rrrr]^f &&&& Y\\
}
$$
is a minimal $\overline \A$-cone decomposition for $f$.  Thus  
each $A_i\in \overline \A$ and each $A_i\maprt{} X_i\maprt{j_i}X_{i+1}$
is a mapping cone sequence.   Since $A_i\in\overline \A$, $\kl_\A(A_i)\leq 1$,
and hence $\LL_\A(j_i) \leq 1$ by Corollary \ref{cor:mappingcones}(2a).
By the Composition Axiom, $\LL_\A(f) \leq n = \LL_{\overline \A}(f)$.
\end{prf}

We next show that the collection $\overline \A$ satisfies the inequality 
of Theorem \ref{thrm:MPT}(1) if and only if $\overline \A$
is closed under both wedges and suspension.
For this it suffices to prove the following corollary.

\bcor
Let $\A$ be any collection and consider commutative diagrams of the form
$$
\xymatrix{
C\ar[d]^{c} && A \ar[ll]_{g} \ar[rr]^{f} \ar[d]^{a} && B\ar[d]^{b}   
\\
C' && A' \ar[ll]_{g'}\ar[rr]^{f'} && B'  . 
\\  }
$$
If the inequality  
\begin{enumerate}
\item $\LL_\A(d) \leq \LL_\A(a) +  \max ( \LL_\A(b),\LL_\A(c))$
\end{enumerate}
of Theorem \ref{thrm:MPT} holds for any such diagram, then $\overline \A$ is 
closed under both wedges and suspension.
\ecor

\begin{prf}
We show that $\overline \A$ is closed under suspension; the proof of the other
assertion is similar.  Let $A\in \overline \A$ and  consider the commutative diagram
$$
\xymatrix{
{*} \ar[d] && A   \ar[ll]\ar[d]\ar[rr] && {*} \ar[d]\\
{*}        && {*} \ar[ll]      \ar[rr] && {*}.\\
}
$$
By (1), we have 
$$
\begin{array}{rclcl}
\kl_\A(\s A) & = & \LL_{\A}(\s A\maprt{} *) \\
& \leq & \LL_{\A}(A\maprt{} *) \\
& = & \kl_\A(A)  & \leq & 1, \\
\end{array}
$$
so $\s A\in \overline \A$ by definition.
\end{prf}

\brmk  {\rm To conclude that $\overline \A$ is closed under suspension,
it suffices to consider only diagrams in which $b = \id $ and $c=\id$,
and to conclude that $\overline \A$ is closed under wedges, 
we only need to consider diagrams with $a=\id$. }
\ermk

\bigskip

\section{Products}

The following is our main result on products of maps.

\bthm\label{thrm:5.1}
Let $\A$ be a collection that is closed under wedges and joins
and let $f:A\maprt{}X$ and $g:B\maprt{}Y$ be maps.  Then

\begin{enumerate}
\item
$\LL_\A(f\cross g) \leq \LL_\A(f) + \LL_\A(g)
+\max( \cl_\A(A) , \cl_\A(B) )$,
\item
$\L_\A(f\cross g) \leq  \L_\A(f) + \L_\A(g)
+\max( \cl_\A(A) , \cl_\A(B) )$.
\end{enumerate}
\ethm

\begin{prf}
In the proof of (1) we write 
$a=\cl_\A(A)$, $b=\cl_\A(B)$, $m=L_\A(f)$ and $n=L_\A(g)$ and assume that $a\geq b$.  

Now consider the $\A$-cone decompositions of $*\maprt{}A$ and $f$
$$
\xymatrix{ 
K_0  \ar[d]           & K_1\ar[d]  &                & K_{a-1}\ar[d]   & K_a\ar[d] && K_{m+a-1}\ar[d]
\\
A_0 \ar[r]\ar@{=}[d]  & A_1\ar[r]  & \cdots\ar[r]   & A_{a-1} \ar[r]  & A_a\dequivarrow \ar[r] & \cdots \ar[r]
& A_{m+a-1}\ar[r] & A_{m+a}\dequivarrow 
\\
{*}\ar[rrrr] &&&& A \ar[rrr]_f &&& X\\
}
$$
and of $*\maprt{}B$ and $g$
$$
\xymatrix{ 
L_0  \ar[d]           & L_1\ar[d]  &                & L_{b-1}\ar[d]   & L_b\ar[d] && L_{m+a-1}\ar[d]
\\
B_0 \ar[r]\ar@{=}[d]  & B_1\ar[r]  & \cdots\ar[r]   & B_{b-1} \ar[r]  & B_b\dequivarrow \ar[r] & \cdots \ar[r]
& B_{n+b-1}\ar[r] & B_{n+b}\dequivarrow 
\\
{*}\ar[rrrr] &&&& B \ar[rrr]_g &&& Y,\\
}
$$
where we identify $A$ with $A_a$ and $B$ with $B_b$.
Since $A_i\sseq A_{i+1}$ and $B_j\sseq B_{j+1}$,
we may define, for $0\leq k \leq n+m+a + b$,  
$$
C_k = A\cross B \ \cup\bigcup_{i+j=k} A_i\cross B_j \ \ \sseq \ \ A_{m+a} \cross B_{n+b} .
$$
Observe that $C_b = A\cross B$, $C_{n+m+a + b} = X\cross Y$
and up to homotopy the composite
$$
C_b \maprt{} C_{b+1} \maprt{} \cdots \maprt{} C_{n+m+a + b}
$$
is $f\cross g$.  From this, we see that it suffices to show that
$\LL_\A(C_k \maprt{} C_{k+1}) \leq 1$ for each $k\geq b$.

For $0\leq i \leq a + m $ and $0\leq j\leq b+ n$, define
$$
P_{ij} = A_i \cross B_j,\ \   T_{ij} = A_i \cross B_{j-1} \cup
A_{i-1}\cross B_j, \ \ {\mathrm{and}}\ \ 
Q_{ij} = C_{i+j-1}\cup P_{ij}.
$$
Then $C_{k+1}$ is obtained as
the pushout of all the maps $C_k \maprt{} Q_{ij}$
with $i+j = k+1$.  
By an induction based on Corollary \ref{cor:pushouts}(4) it follows that
$$ \LL_\A(C_k\maprt{}C_{k+1})\leq
\max( \LL_\A(C_k \maprt{} Q_{ij})).
$$
Thus, it suffices to show that for $i+j = k+1$,
$\LL_\A(C_k \maprt{} Q_{ij})\leq 1$.  Applying  
Corollary \ref{cor:pushouts}(1) to the pushout
diagram
$$
\xymatrix{
T_{ij}\ar[rr]\ar[d] && C_k \ar[d]
\\
P_{ij}\ar[rr] && Q_{ij},
}
$$  
we have
$\LL_\A(C_k \maprt{} Q_{ij})\leq \LL_\A( T_{ij}\maprt{}P_{ij})$.
According to a result of Baues \cite{Baues1} (see also 
\cite{Stanley2}),
there is a mapping cone sequence
$$
K_{i-1} \ast L_{j-1} \maprt{} T_{ij}\maprt{} P_{ij} 
$$
when $i,j>0$,   a mapping cone sequence
$$ \xymatrix{
L_{j-1} \ar[r] & T_{0j} \ar@{=}[d]\ar[r] & P_{0j}\ar@{=}[d]  
\\
                & B_{j-1}                  & B_j 
\\ }
$$
when $i=0$ and a mapping cone sequence
$$ \xymatrix{
K_{i-1} \ar[r] & T_{i0} \ar@{=}[d]\ar[r] & P_{i0}\ar@{=}[d] 
\\
  & A_{i-1}                  & A_i 
\\ }
$$
when $j=0$.
Since $\A$ is closed under joins,  
$\LL_\A(T_{ij}\maprt{} P_{ij})\leq 1$, and this completes the proof of (1).

For (2) we take $f'$ to be a map which dominates $f$, 
has the same domain and such that $\L_\A(f) = \LL_\A(f')$, and $g'$ 
is similarly chosen for $g$ (Lemma \ref{lem:strict}).  Then (2) is a
consequence of (1) since $f'\cross g'$ dominates $f\cross g$.
\end{prf}

\bcor\label{cor:5.2}
If $\A$ is closed under wedges and joins, then
\begin{enumerate}
\item \begin{enumerate}
\item  $\cl_\A (X\cross Y) \leq \cl_\A(X) + \cl_\A(Y)$,
\item  $\kl_\A(X \cross Y) \le  \kl_\A(X) + \kl_\A(Y) + \max(\cl_\A(X),\cl_\A(Y))$;
\end{enumerate}
\item \begin{enumerate}
\item $\cat_\A(X\cross Y) \le \cat_\A(X) + \cat_\A(Y)$,
\item $\kit_\A(X\cross Y) \le \kit_\A(X) + \kit_\A(Y) + \max(\cl_\A(X),\cl_\A(Y))$.
\end{enumerate}\end{enumerate}
\ecor

\brmk
{\em In the case $\A=\{ \mathrm{ all\ spaces} \}$,  Corollary \ref{cor:5.2}(2a)
is a classical result due to Bassi \cite{Bassi}.
Part (1a) has been obtained by Takens \cite{Takens},
Clapp and Puppe \cite{ClP2}, and Cornea \cite{Cornea5}.}
\ermk

It is possible to improve the inequalities in Corollary \ref{cor:5.2} by imposing
stronger conditions on the collection $\A$.  To illustrate this, we state
and sketch a proof of Proposition \ref{prop:below} below.  We say that a collection
$\A$ is a \term{$\smsh$-ideal} if for any $A\in\A$ and any space $B$,
the smash product $A\smsh B\in \A$.

\bprp\label{prop:below}
If $\A$ is a $\smsh$-ideal and is closed under wedges and suspensions, then
\begin{enumerate}
\item  $\kl_\A(X\times Y)\leq \kl_\A(X) + \kl_\A(Y)$ and 
\item  $\kit_\A(X\times Y)\leq \kit_\A(X) + \kit_\A(Y)$.
\end{enumerate}
\eprp

\begin{prf}
We only prove (1) since the proof of (2) is similar.  
By applying Corollary \ref{cor:mappingcones}(4) 
to the sequence $X\wdg Y\maprt{}X\times
Y\maprt{}X\smsh Y$ we conclude that 
$\kl_\A(X\times Y)\leq \kl_\A(X\wdg Y) +  \kl_\A(X\smsh Y)$.
By Corollary \ref{cor:susppushouts}(2), 
$\kl_\A(X\wdg Y)\leq \max(\kl_\A(X),\kl_\A(Y))$.  
Furthermore, a simple argument using the
fact that $\A$ is a $\smsh$-ideal shows that $\kl_\A(X\smsh Y)\leq
\min(\kl_\A(X),\kl_\A(Y))$.  This completes the sketch of the proof.
\end{prf}

\bigskip

\section{Pullbacks and Fibrations}

We prove a result on pullbacks which yields inequalities for the 
$\A$-cone length and $\A$-category of the spaces  which appear in a fiber sequence.  

We begin with a lemma which may be of independent interest.  
In the proof we denote the half-smash $(X \cross Y)/X$ by $X \SDP Y$ 
and the quotient map by $q:X \cross Y \maprt{} X\SDP Y$. 

\blem\label{lem:halfsmash} {\rm (Cf. \cite[Ex.\thinspace 5.4]{Marcum})}
Let $\A$ be a collection which is closed under joins and let $A\in \A$.  
If $p_2:A\cross B\maprt{}B$ is the projection, then
\begin{enumerate}
\item[(1)] $\LL_\A(p_2) \leq \cl_\A(B) + 1$
\item[(2)] $\L_\A(p_2) \leq \cat_\A(B) + 1$.
\end{enumerate}
\elem

\begin{prf}
Consider the map $p:A\SDP B\maprt{}B$ induced by $p_2$.   
The main step in the proof is to show 
$$ \LL_\A(p) \leq \cl_\A(B)\qquad {\mathrm{and}}\qquad   \L_\A(p) \leq \cat_\A(B).
$$
Suppose we have a diagram:
$$
\xymatrix{ 
L_0  \ar[d]   & L_1\ar[d] & &L_{n-1}\ar[d]
\\
{*}=B_0 \ar[r]^-{j_0} & B_1\ar[r]^-{j_1} & \cdots\ar[r]   & B_{n-1} \ar[r]^-{j_{n-1}} & B_n \\
}
$$
and a map $f_n:B_n\maprt{}B$ with $\LL_i \in A$.
Define $D_i$ as the homotopy pushout in the diagram
$$
\xymatrix{      
A\SDP B_i\ar[rr]^{q_i}\ar[d]_{\id\sdp  f_nj_{n-1}...j_i }  && B_i\ar[d]^{r_i} 
\\
A\SDP B\ar[rr]^{s_i} && D_i, \\
}
$$
where $q_i$ is the projection.
Then there are maps $k_i:D_i\maprt{}D_{i+1}$ with $k_i s_i =s_{i+1}$.
When $i = 0$ we have $D_0 = A\SDP B$ and when $i =n$ we have 
$$
\xymatrix{
A\SDP B_n \ar[rr]^{q_n} \ar[d]_{\id\sdp f_n}&& B_n \ar[d]^{r_n}\\
A\SDP B   \ar[rr]^{s_n}                        && D_n. \\
}
$$
From the above diagram and the maps 
$$
 f_n:B_n\maprt{}B  \quad  {\mathrm{and}} \quad  p:A\SDP B\maprt{}B , 
$$
we obtain a map $g_n:D_n\maprt{}B$ such that $g_ns_n = p$ and $g_nr_n = f_n$.  
It then follows that
$$
g_nk_{n-1}\cdots k_0 = p.
$$
Now we prove (1).  Suppose $f_n$ is a homotopy equivalence 
so our given decomposition is an $\A$-cone
decomposition of $B$ of length $n$.  By the previous 
homotopy pushout diagram, $r_n$ is a homotopy equivalence 
and from $g_nr_n = f_n$ we obtain that $g_n$ is a 
homotopy equivalence.   Since $g_nk_{n-1}\ldots k_0 = p$, we get $\LL_A(p) \leq \LL_\A(k_0) +
\cdots +\LL_\A(k_{n-1})$.  To  complete the proof that $\LL_\A(p)\leq n =\cl_\A(B)$, 
it suffices to show  that $\LL_\A(k_i) \leq 1$.  But $L_i\maprt{}B_i\maprt{j_i}B_{i+1}$
is a mapping cone sequence and so 
$$
A\SDP \LL_i\maprt{}A\SDP B_i\maprt{}A\SDP B_{i+1}
$$ 
is also a mapping cone sequence.  Thus we have a commutative  diagram
$$
\xymatrix{
{*}\ar[d]&& A\SDP L_i\ar[ll]\ar[rr]\ar[d] && \LL_i\ar[d]
\\
A\SDP B\ar[d]&& A\SDP B_i\ar[ll]\ar[rr]^{q_i}\ar[d]&& B_i\ar[d]
\\
A\SDP B&& A\SDP B_{i+1}\ar[ll]\ar[rr]^{q_{i+1}} && B_{i+1}.\\
}
$$
Since each column is a cofiber sequence, $P\maprt{}D_i\maprt{k_i}D_{i+1}$ is a cofiber 
sequence by the Four Cofibrations Theorem, where $P$ is the homotopy pushout of the top
line.  However it is easily seen that  $P = A*L_i$, the join of $A$ and $\LL_i$.  But $P\in \A$
since $\A$ is closed under joins.  Thus $\LL_\A(k_i)\leq 1$ and so $\LL_\A(p)\leq n =\cl_\A(B)$.  
Part 1 of the lemma now follows by factoring $p_2:A\cross B\maprt{}B$ as
$$
A\cross B\maprt{q}A\SDP B\maprt{p}B.
$$
Since $A\maprt{}A\cross
B\maprt{q}A\SDP B$ is mapping cone sequence with $A \in \A$,
$\LL_\A(q)\le 1$.  Thus $\LL_\A(p_2)\leq \LL_\A(p)+\LL_\A(q)\leq \cl_\A(B) +1$
by the Composition Axiom.

The proof of (2) is similar.  Instead of taking an $\A$-cone decomposition
of $B$, we take an $\A$-category decomposition of $B$ of length $n$.  
Thus instead of having
$f_n:B_n\maprt{}B$ a homotopy equivalence, we have a map $s: B\maprt{}B_n$ with 
$f_n s \simeq\id$.   We define $\sigma:B\maprt{}D_n$ by $\sigma = r_n s$.  Then 
the following are easily checked:
$$
\begin{array}{rlrlrl}
{\mathrm (a)} & g_n \sigma \simeq \id  \qquad               & 
{\mathrm (b)} & \sigma p \simeq  k_{n-1}\cdots k_0  \qquad   &
{\mathrm (c)} & g_n k_{n-1}\cdots k_0 \simeq p  .     \\
\end{array}
$$
Using the maps $(\id, \sigma)$ and $(\id,g_n)$ we see that $p$ is homotopy dominated by
$k_{n-1}\cdots k_0$. Therefore $\L_\A(p) \leq  \L_\A(k_0) + \cdots +\L_\A(k_{n-1})$.  
The rest of the proof is the same as the proof of (1), 
using $\L$ for $\LL$.    
\end{prf}

Now we prove our pullback theorem.

\bthm\label{thrm:5.3}
Let $\A$ be a collection that is closed under wedges and joins and let
$$
\xymatrix{
A\ar[rr]\ar[d] && B\ar[d]\\
C\ar[rr] && D\\  }
$$
be a pullback diagram.  Let $B\maprt{}D$ be a fibration with fiber
$F$.  Then
\begin{enumerate}
\item $\LL_\A(A\maprt{} B) \leq \LL_{\A} (C\maprt{}D) ( \cl_\A (F) +1)$;
\item $\L_\A(A\maprt{} B) \leq
\L_{\A} (C\maprt{}D) ( \cat_\A (F) +1)$.
\end{enumerate}
\ethm

\begin{prf}
We prove (1).  Let
$$
\xymatrix{ 
K_0  \ar[d]           & K_1\ar[d]  &                & K_{n-1}\ar[d]
\\
C_0 \ar[r]\ar@{=}[d]  & C_1\ar[r]  & \cdots\ar[r]   & C_{n-1} \ar[r]  & C_n\dequivarrow 
\\
C\ar[rrrr] &&&& D\\
}
$$
be a minimal $\A$-cone decomposition for $C\maprt{}D$.
For $0\leq i\leq n$, define $B_i$
to be the pullback indicated by the square
$$
\xymatrix{
B_i\ar[rr]\ar[d] &&  B\ar[d] \\
C_i\ar[rr] && D.\\  }
$$
Thus $B_0 = A$, $B_n\equiv B$ and we obtain maps $B_i\maprt{}B_{i+1}$.
With these identifications, the composition
$B_0\maprt{}B_1 \maprt{}\cdots \maprt{} B_n$ is simply $A\maprt{}B$.
Hence, it suffices to show that $\LL(B_i\maprt{} B_{i+1})\leq
\cl_{\A}(F) +1$.  Consider the cube diagram
$$
\xymatrix{
& K_i \cross F \ar[rr]\ar[ld]\ar'[d][dd]      && F \ar[dd] \ar[ld]  
\\
B_i \ar[rr] \ar[dd]                  && B_{i+1} \ar[dd] 
\\
& K_i\ar'[r][rr]\ar[ld]             && {*} \ar[ld]
\\
C_i\ar[rr]                          && C_{i+1}.
}
$$
In this diagram, the bottom square is a homotopy pushout and
the sides are pullbacks.  This assertion is obvious for all squares except the
left side square
$$
\xymatrix{
K_i\cross F \ar[rr]\ar[d]    && B_i \ar[d] \\
K_i \ar[rr] && C_i . \\ }
$$   
To see that this  is a pullback square, let $P$ be the 
pullback  of $B_i\maprt{}C_i\maplft{}K_i$.  Then $P$ is
also the pullback of $B_{i+1}\maprt{}C_{i+1}\maplft{}C_i\maplft{}K_i$.  
Since the composite $C_{i+1}\maplft{}C_i\maplft{}K_i$ is the constant map, the  latter 
pullback is $K_i \cross F$.

Now by the Mather's second cube
theorem \cite{Mather}, the top square is a homotopy pushout.
By  Corollary \ref{cor:pushouts}(1),
$\LL_\A(B_i\maprt{} B_{i+1}) \leq \LL_\A ( K_i\cross F \maprt{} F )$.
By Lemma \ref{lem:halfsmash}(1), $\LL_\A(K_i \cross F\maprt{} F)  \leq \cl_\A(F) + 1$.  
This proves (1).

The proof of (2) is similar and uses Lemma 6.1(2) and we omit it.
\end{prf}

\bcor\label{cor:5.4}
Let $\A$ be a collection that is closed under wedges and joins
and let $F\maprt{} E\maprt{}B$ be a fibration.  Then
\begin{enumerate} 
\item $\cl_\A(E) + 1 \leq (\cl_{\A} (B) + 1)(\cl_\A(F) +1)$;
\item $\cat_\A(E) + 1 \leq (\cat_{\A} (B) + 1)(\cat_\A(F) +1)$.
\end{enumerate}
\ecor

\begin{prf}
We prove (1).  Applying Theorem \ref{thrm:5.3} to the  pullback square
$$
\xymatrix{
F\ar[rr] \ar[d] && E\ar[d]  \\
{*} \ar[rr]     && B,        \\    }
$$
we obtain $\LL_\A(F\maprt{} E) \leq  \cl_\A(B) (\cl_\A(F) + 1)$.
Now the Composition Axiom shows that 
$$
\cl_\A(E) \leq \cl_\A(F) + \LL_\A(F\maprt{}E),
$$
so
$$
\begin{array}{rcl}
\cl_\A(E) + 1                               & \leq &
\cl_\A(F) + \LL_\A(F\maprt{}E) + 1  \\      & \leq &
\cl_\A(B) (\cl_\A(F) + 1) + (\cl_\A(F) + 1) \\ & = &
(\cl_{\A} (B) + 1)(\cl_\A(F) +1).
\end{array}
$$
\end{prf}

\brmk {\rm
In the special case $\A=\{ \mathrm{all\ spaces}\}$
we retrieve Varadarajan's result \cite{V}  
$$
\cat(E) + 1 \leq (\cat(B) +1) (\cat(F) +1).
$$
Hardie has obtained a further improvement in \cite{Hardie2}, but that 
involves a different notion of the category of a map 
from the one we consider here \cite{B-G,Fox}.  }
\ermk

\bigskip

\section{Miscellaneous Results and Problems}

In this section we consider several topics.  
We first establish some elementary, but useful, facts about 
$\LL_\A$ and $\L_\A$.  We then show
that some known results for the collection 
$\A = \{{\mathrm{all\ spaces}}\}$ do not hold for an
arbitrary collection $\A$.  
Finally, we conclude the section by stating a number of open
questions and discussing them briefly.

We begin with a few elementary results.

\bprp  Let $f:X\maprt{}Y$ and let $\A$ be any 
collection.  Then
\begin{enumerate}
\item $\LL_\A(f) = 0$  if and only if $f$ is a homotopy equivalence;
\item $\L_\A(f) = 0$ if and only if $f$ is a homotopy equivalence.
\end{enumerate}
\eprp 

\begin{prf}  We   prove (1) and (2) at the same time.  By the axioms, if $f$ is a homotopy
equivalence, then $\L_\A(f) = \LL_\A(f)  = 0$.  Conversely, define a
function $\l_\A$ by 
$$
\l_\A (f) = \Biggl\{  
\begin{array}{ll}
0 & {\mathrm {if}}\ f\ {\mathrm{is\ a\ homotopy\ equivalence}}\\
1 & {\mathrm {otherwise}}.\\
\end{array}
$$
It is trivial to check that $\l_\A$ satisfies the $\A$-category
axioms, so 
$$
\l_\A(f)   \leq \L_\A(f)\leq \LL_\A(f)
$$ 
for every map $f$.  Consequently, if $f$ is not a homotopy equivalence, 
then $\LL_\A(f) \geq\L_\A(f) \geq \l_\A(f) =1$.
\end{prf}

\bprp\label{prop:wdg}
Let $f:X\maprt{}Y$ and $g:X'\maprt{}Y'$ be maps and let 
$\A$ be a collection that is closed under wedges.  Then
\begin{enumerate}
\item[(a)]  $\LL_\A(f\wdg g)  \leq \max(\LL_\A(f), \LL_\A(g))$;
\item[(b)]  $\L_\A(f\wdg g) = \max(\L_\A(f), \L_\A(g))$.
\end{enumerate}
\eprp

\begin{prf}
Since both $X \wdg X'$ and $Y\wdg  Y'$ are homotopy pushouts, the inequality
$\L_\A(f \wdg g) \leq \max(\L_\A(f), \L_\A(g))$ is a consequence of Theorem 
\ref{thrm:MPT}.  This same argument shows  $\LL_\A(f \wdg g) \leq \max(\LL_A(f), \LL_A(g))$.  
The reverse inequality for $\L_\A(f \wdg g)$ follows since $f$ and $g$ are both retracts
of $f \wdg g$.
\end{prf}

An example due to Dupont \cite{dupont} shows that 
equality does not generally hold in (a).  Other examples
can be found in \cite{Stanley}, where spaces $X_n$ with 
category $n$ and cone length $n+1$ are constructed.  
According to an observation of Ganea \cite{Takens}
(see also \cite{Cornea2}), this implies that 
there is a space  $A$ such that $\cl( X_n\wdg \s A ) = \cat(X_n)$.
If we let $\A=\{ {\mathrm{all\ spaces}}\}$,
$f:*\maprt{} X_n$ and $g:*\maprt{}\s A$,
then we have
$$
\LL_\A( f ) = \cl (X_n) > \cat(X_n) = \cl (X_n\wdg \s A) = \LL_\A(f\wdg g).
$$

\bcor  Let $X$ and $Y$ be spaces and $\A$ a collection that is
closed under wedges.  Then 
\begin{enumerate}
\item[(a)]  $\LL_\A(X\maprt{*}Y) \leq \max ( \kl_\A(X), \cl_\A (Y))$;
\item[(b)]  $\L_\A(X\maprt{*}Y) = \max ( \kit_\A(X), \cat_\A (Y))$.
\end{enumerate}
\ecor

\begin{prf}
The trivial map $X\maprt{*}Y$ is the wedge of the maps $X \maprt{} *$
and $*\maprt{}Y$.  The result follows from Proposition \ref{prop:wdg}.
\end{prf}

Next we turn to some known results for the collection 
$\A = \{{\mathrm{all\ spaces}}\}$.   For this collection we 
delete the subscript $\A$ and
write $\LL_\A$ as $\LL$, $\L_\A$ as $\L$, etc.  

For any map $f:X\maprt{}Y$, it has been proved that $\LL(f) \leq \cl(Y) + 1$ 
\cite{Marcum} and $\L(f) \leq \cat(Y) + 1$ \cite{Cornea2}.  
We show that this may not be true for an arbitrary collection $\A$.

\bexm {\rm  By Corollary \ref{cor:lowerbound},
$\kl_\A(X) \leq \LL_\A(f) + \kl_\A(Y)$.
Thus if $L_\A (f) \le \cl_\A(Y) + 1$ were true, we would have
$$
\kl_\A( X )\leq  \cl_\A(Y) + \kl_\A(Y) + 1,
$$
for {\it any} $X$ and $Y$.  This cannot hold for any collection 
$\A$ such that there are spaces $X$ with arbitrarily large 
killing length (e.g., for $\A = \S$ or $\s$).
The analogous observation holds for $\L_\A$.}
\eexm

Another classical result concerns the homotopy pushout square
$$
\xymatrix{
A\ar[rr] \ar[d] && B\ar[d]\\
C\ar[rr] && D .\\ }
$$
It has been  shown that $\cl(D) \leq \cl(B) + \cl(C) + 1$  \cite{Hardie}.
We show that this is not true for $\A = \S$, the collection of wedges of spheres.

\bexm {\rm
Consider the homotopy pushout
$$
\xymatrix{
\cp^{t}\ar[rr] \ar[d] && {*}\ar[d]\\
{*}\ar[rr] && \s\cp^{t}.\\ }
$$   
As $t$ increases, the length of the longest nontrivial composition  of 
Steenrod squares in $H^*(\s \cp^t; \ZZ_2)$ also becomes arbitrarily large.  
It follows from \cite[Prop.\thinspace 7.5]{ASS}
that $\cl_\S(\s \cp^t)$ increases as $t$ increases.  
This contradicts the $\S$-analog of Hardie's result.  }
\eexm

We conclude the paper by stating and discussing three open problems.

\medskip

\noindent
{\bf Problem 7.6}\quad   
We have seen in \cite[Prop.\thinspace 7.3]{ASS} that for certain collections $\A$, 
$\wcat(X) \leq 2^{\kl_\A(X)} - 1$,
where $\wcat(X)$ is the weak category of X (see \cite{Gilbert,James}).  
Since $\wcat(X) \leq\cat(X) \leq \cat_\A(X)$ for any collection $\A$, 
it is reasonable to ask for which collections $\A$ is 
$\cat_\A(X) \leq 2^{\kl_\A(X)} - 1$.
Of course $\A$ must not be $\{{\mathrm{all\ spaces}}\}$,
since   $\kl_\A(X) \leq 1$ for every space $X$ in that case.  We note
that the conjecture has been verified in the case $\A = \S$ and 
$X = S^n_1 \cross\cdots \cross S^n_r$ \cite[Prop.\thinspace 6.2]{AMS}.  
Other evidence for the conjecture in the case $\A = \Sigma$,  
the collection of suspensions, has been given in
\cite{AS}, where a weaker form of  this problem has been posed 
\cite[\sec 7,\thinspace No. 5]{AS}.

\medskip

\noindent
{\bf Problem 7.7}\quad   
Given $f:X\maprt{}Y$.  Is $\LL_\A(f) \leq \kl_\A(X) + \cl_\A(Y)$, 
and is $\L_A(f) \leq  \kit_\A(X) + \cat_\A(Y)$?

We discuss the evidence in the case of $\LL_\A$ (the discussion
is analogous for $\L_\A$).
First of all, if $C$ is the cofiber of $f$ it is true that 
$\cl_\A(C) \leq \kl_\A(X) + \cl_\A(Y)$ 
(Corollary \ref{cor:mappingcones}(3a)) 
and also $\cl_\A(C) \leq \LL_\A(f)$
(Corollary \ref{cor:mappingcones}(1a)). 
Secondly, we have that  $\LL_\A(f) \leq \cl_\A(X) + \cl_\A(Y)$ 
(Corollary \ref{cor:4.10}(1a)) and $\kl_\A(X)\leq \cl_\A(X)$. 
Finally, when $\A=\{ \mathrm{all\ spaces}\}$ then $\kl_\A(X)=1$ for 
every $X$, and in this case it is known that $L(f) \leq \cl(Y) + 1$ 
\cite{Marcum}.

\medskip

\noindent
{\bf Problem 7.8}\quad  
It is well known that $\cl(X)\leq \cat(X)+1$ \cite{Takens}.
If $\A$ is a collection different from $\{ \mathrm{all\ spaces}\}$,
is there an upper bound for $\cl_\A(X)$ in terms of 
$\cat_\A(X)$? This question was asked by Scheerer-Tanr\'e in \cite{S-T}.
Analogously, is there an upper bound for $\kl_\A$ in terms of $\kit_\A$?

We can show that $\kit_\s (X) \leq 1$ implies
$\kl_\s (X) \leq 3$ as follows.  If $\kit_\s(X) \leq 1$
then there is a mapping cone sequence $A\maprt{} X \maprt{*} Y$
with $A\in\s$.  It follows that $\s A = Y \wdg \s X$,
and so there is a retraction
map $\alpha : \s A \maprt{} Y$.
The cofiber of $\alpha $ is $\s^2 X$, and hence we have a 
decomposition
$$
\xymatrix{
A\ar[d] &&  \s A \ar[d] && \s^2 X \ar[d]\\
X\ar[rr] && Y \ar[rr] && \s^2 X \ar[rr] && {*}.\\ }
$$
This proves that $\kl_\s(X) \leq 3$.

%
%

Department of Mathematics

Dartmouth College

Hanover, NH 03755 

U.~S.~A

Martin.Arkowitz@Dartmouth.edu

\bigskip

University of Alberta, 

Edmonton, AB  T6G 2G1

CANADA

dstanley@sirius1.math.ualberta.ca

\bigskip

Department of Mathematics

Western Michigan University

1903 W. Michigan Ave.

Kalamazoo, MI 49008

U.~S.~A

Jeffrey.Strom@wmich.edu

\end{changemargin}

\end{document}